\documentclass[11pt]{amsart}
\usepackage{hyperref}
\newtheorem{thm}{Theorem}[section]

\theoremstyle{definition}

\theoremstyle{remark}

\numberwithin{equation}{section}
\newcommand{\Real}{\mathbb R}

\begin{document}

\title{Scattering theory for p-forms on hyperbolic real space}
\author{Francesca Antoci}
\address{Dipartimento di Matematica, Politecnico di Torino}
\email{antoci@dimat.unipv.it}

%\thanks{}

%\date{}%
%\dedicatory{}%
%\commby{}%
% ----------------------------------------------------------------
\begin{abstract}
Due to spectral obstructions, a scattering theory in the Lax-Phillips sense for the wave equation
for differential $p$-forms on $H^{n+1}$ cannot be developed. As a consequence, Huygens' principle
for the wave equation in this context does not hold. \par
If we restrict the class of forms and we consider the case of coclosed $p$-forms on $H^{n+1}$,
 when $n=2p$, Huygens' principle does hold and thus in this case incoming and outgoing subspaces
 can be constructed.
\end{abstract}
\maketitle
% ----------------------------------------------------------------
\section{Introduction}
Scattering theory for the wave equation has been developed by Lax and Phillips essentially for
functions, at first in euclidean setting (see \cite{Lax-Phillips} and the references therein), and
more recently, after the papers of Faddeev and Pavlov (see \cite{Faddeev}) and the construction of
a Fourier analysis for functions on symmetric spaces due to Helgason (see \cite{Helgason}), in
non-euclidean setting. Particular attention has been devoted in \cite{Lax-Phillips2},
\cite{Lax-Phillips3}, \cite{Lax-Phillips4}, \cite{Lax-Phillips5} to the wave equation for
functions on the hyperbolic space $H^{n+1}$ and on quotients of $H^{n+1}$ by groups of
isometries.\par
The purpose of this paper is to extend the Lax-Phillips theory to the case of the wave equation
for $p$-forms, or more generally for sections of vector bundles, on complete, non compact
Riemannian manifolds. The basic motivation is connected, on one hand, with the spectral properties
of the Laplace-Beltrami operator on differential forms, and hence in some way with the ``geometry''
of the manifold, and on the other side with the validity of Huygens' principle for the wave
equation, and therefore with the ``analysis'' that can be developed on the manifold.\par
A general approach to these problems on a general manifold turns out to be extremely difficult
 to deal with, for the lack of information about the behaviour of Laplace-Beltrami operator on
 a generic complete non compact Riemannian manifold, and of analytic tools mimicking the Radon
 and Fourier transforms.\par
These facts motivate our choice to investigate the problem at first on the hyperbolic space
$H^{n+1}$. In this case, the spectrum of the Laplace-Beltrami operator on $p$-forms on $H^{n+1}$
has been explicitely computed by Donnelly (see \cite{Donnelly}). Moreover in recent years
(see \cite{Branson-Olafsson-Schlichtkrull}) it has been developed a Radon transform for $p$-forms
 on $H^{n+1}$, providing an essential tool to the study of the wave equation.\par
In section 2, we recall some essential facts about $H^{n+1}$ and its spectral properties.
In section 3, we prove that it is not possible to develop a scattering theory in the Lax-Phillips
sense for the wave equation for all $p$-forms on $H^{n+1}$. In section 4, we show that, as a
consequence, Huygens' principle for the wave equation in this context does not hold. In section
5 we show that since  Huygens' principle does hold for the wave equation for coclosed $p$-forms
when $n=2p$, a scattering theory in the Lax-Phillips sense can be developed if we impose some
restrictions on the class of forms considered. This result implies some conclusions about the
spectrum of the Laplace-Beltrami operator on coclosed $p$-forms on $H^{n+1}$, when $n=2p$. \par
This paper summarizes some results of my PhD Thesis, compiled under the direction of Prof.
Edoardo Vesentini, to whom I am deeply grateful for his constant and careful support.

\section{Preliminaries}

We will denote by  $H^{n+1}$ the real hyperbolic space of dimension $n+1$, either in the
half-space model
$$ H^{n+1}= \left\{ (x_1,...,x_{n+1}) \in R^{n+1}\; \mid \;x_{n+1}>0\,\right\}$$
endowed with the Riemannian metric
$$ds^2= \frac{dx_1^2+dx_2^2+...+dx_{n+1}^2}{x_{n+1}^2},$$
or in the Poincar\'e model
$$H^{n+1}= \left\{ (x_1,...,x_{n+1}) \in {\Real}^{n+1}\; \mid \; |x|^2 <1 \,\right\}$$
with the Riemannian metric
$$ d\sigma^2(x)= \left( \frac{2}{1- |x |^2} \right)^2 (dx_1^2+...+dx_{n+1}^2),$$
where $|x|^2=x_1^2+x_2^2+...+x_{n+1}^2$.\par \bigskip
We will denote by $\Lambda^p(H^{n+1})$ the space of smooth $p$-forms on $H^{n+1}$, and by
$\Lambda_c^p(H^{n+1})$ the space of smooth compactly supported $p$-forms on $H^{n+1}$.
As usual,
$L^2_p(H^{n+1})$ will be the Hilbert space of square integrable $p$-forms on $H^{n+1}$, that is
the completion  of $\Lambda_c^p(H^{n+1})$ with respect to the $L^2$ norm induced by the
Riemannian metric:
$$\|\omega\|_{L^2}^2= \int_{H^{n+1}} \omega \wedge * \omega.$$

Since $H^{n+1}$ is a complete Riemannian manifold, the Laplace-Beltrami operator $\Delta_p : \Lambda^p_c(H^{n+1}) \rightarrow \Lambda^p_c(H^{n+1})$ is an essentially selfadjoint operator in $L^2_p(H^{n+1})$. We will denote its closure again by the same symbol $\Delta_p$. \par
H. Donnelly in \cite{Donnelly} has explicitely computed the spectrum of the Laplace-Beltrami operator on $p$-forms on $H^{n+1}$. He has proved that $\Delta_p$ has point spectrum if and only if $n+1= 2p$. In this case, the only eigenvalue of $\Delta_p$ is zero, occurring with infinite multiplicity, and the continuous part of the operator is unitarily equivalent to the operator
$$ L^2(\Real^+,dr,  M \oplus  M) \longrightarrow L^2({\Real}^+,dr,  M
\oplus  M)$$
$$f_1(x) \oplus f_2(x) \longmapsto\left( \frac{1}{4}+ x^2\right) f_1(x)\oplus \left( \frac{1}{4}+ x^2 \right) f_2(x) ,$$
where $M$ is a Hilbert space of countably infinite dimension and $dr$ is the Lebesgue measure; hence, if $n+1= 2p$, the absolutely continuous spectrum of $\Delta_p$  is $\sigma_{ac}(\Delta_p)= [\frac{1}{4}, +\infty)$, with constant multiplicity.\par
\noindent If $p \not=  \frac{n+1}{2}$, the point spectrum is empty, and $\Delta_p$ is unitarily equivalent to the operator
$$ L^2({\Real}^+,dr,  M \oplus  M) \longrightarrow L^2({\Real}^+,dr,  M\oplus  M)$$
$$f_1(x) \oplus f_2(x) \longmapsto \left( (\frac{n-2p}{2})^2 + x^2 \right)f_1(x) \oplus \left( (\frac{n-2p+2}{2})^2 + x^2 \right) f_2(x).$$
Hence $\sigma(\Delta_p)= \sigma_{ac}(\Delta_p)= [\mbox{min}\{(\frac{n-2p}{2})^2,(\frac{n-2p+2}{2})^2\}, +\infty)$, and the multiplicity of the spectrum is not constant.\par
\bigskip
\noindent

In the following we will suppose $p \not=  \frac{n+1}{2}$. \par

\section{General case}

Consider the Cauchy problem for the wave equation for $p$-forms on $H^{n+1}$:
\begin{equation}\label{Cauchy} \left\{ \begin{array}{l} \Delta_p \omega +
\partial^2_{tt}\omega=0 \\
                    \omega(0)= \omega_1 \\
                    \partial_t \omega(0)= \omega_2, \end{array} \right.\end{equation}
with $\omega_1, \omega_2$ smooth and compactly supported.
For every $t\geq 0$
\begin{equation}\label{energy}<\Delta_p \omega(t), \omega(t)>+ <\partial_t \omega(t),
\partial_t \omega(t)> = <\Delta_p \omega_1, \omega_1>+ <\omega_2, \omega_2> ,\end{equation}
because
$$\partial_t \left( <\Delta_p \omega(t), \omega(t)>+ <\partial_t \omega(t),
\partial_t \omega(t)> \right) =$$
$$=<\Delta_p (\partial_t \omega),\omega> + <\Delta_p \omega, \partial_t \omega> +$$
$$+ <\partial_{tt}^2 \omega, \partial_t \omega>+
<\partial_t \omega,\partial^2_{tt} \omega>=0.$$
Hence, the energy preserved by the system is given by
$$\|(\omega_1, \omega_2)\|^2_{ H} = <\Delta_p \omega_1, \omega_1>+ <\omega_2, \omega_2> \,\geq 0,$$
which defines indeed a norm on $\Lambda^p_c(H^{n+1})$, because $\|(\omega_1, \omega_2)\|_{ H}=0$
implies:
$$\left\{ \begin{array}{l} d\omega_1=0 \\
                    \delta \omega_1=0 \\
                    \omega_2=0,\\ \end{array} \right.  $$
that is
$$ \left\{ \begin{array}{l} \Delta_p \omega_1=0 \\
                    \omega_2=0. \\ \end{array} \right.$$
Since $p \not= \frac{n+1}{2}$ and therefore $0$ is not in the point spectrum, $\Delta_p
\omega_1 =0$ implies $\omega_1=0$, showing in conclusion that

$$\|(\omega_1, \omega_2) \|_{ H}=0 \; \Longrightarrow \; (\omega_1, \omega_2)=(0,0).$$
\par \bigskip
It is therefore meaningful to take the completion $H$ of the space of smooth compactly supported
initial data with respect to this latter norm, which is induced by the scalar product
\begin{equation}\label{scalarproduct}<(\omega_1,\omega_2), (\tilde{\omega}_1, \tilde{\omega}_2)>=
<\Delta_p\omega_1, \tilde{\omega}_1>+<\omega_2, \tilde{\omega}_2>.\end{equation}
The operator $A$ defined on
$$D(A)= \Lambda^p_c(H^{n+1}) \oplus \Lambda^p_c(H^{n+1}) $$
by
\begin{equation}\label{A} A= \left( \begin{array}{cc} 0 &I \\
                               -\Delta_p &0 \end{array} \right) \end{equation}
is essentially skewselfadjoint in $H$. Hence its closure $\overline{A}$ is the generator of a
one-parameter group of unitary operators $U(t)$
on $H$ which maps every initial data into the value at time $t$ of the corresponding solution
of the wave equation for $p$-forms on $H^{n+1}$.\par
The operator $i\overline{A}$ is essentially selfadjoint and can be shown to be unitarily
equivalent to the operator
$$\left( \begin{array}{cc} \sqrt{\Delta_p} &0 \\
                             0 &- \sqrt{\Delta_p} \end{array} \right)$$
on $L^2_p(H^{n+1}) \oplus L^2_p(H^{n+1})$ (see \cite{Reed-Simon}).\par
Now, thanks to Donnelly's result,  $\sqrt{\Delta_p}$ is unitarily equivalent to the direct
sum of multiplication operators
$$L^2({\Real}^+,dr, M\oplus M) \longrightarrow L^2({\Real}^+,dr, M \oplus M)$$
$$f_1(x) \oplus f_2(x) \mapsto \left( \sqrt{(\frac{n-2p}{2})^2 + x^2}\right)f_1(x) \oplus
\left(\sqrt{(\frac{n-2p+2}{2})^2 + x^2}\right)f_2(x).$$

As a consequence,  if $n+1 \not= 2p$ the spectrum of $i\overline{A}$ is purely absolutely
continuous, equal to
$$\left( - \infty, \max \left\{ - \left|\frac{n-2p}{2}\right|, - \left| \frac{n-2p+2}{2}\right|
\right\}\right]$$ $$\bigcup$$ $$\left[\min\left\{ \left|\frac{n-2p}{2}\right|,\left|
\frac{n-2p+2}{2}\right| \right\}, +\infty\right)$$
and does not have constant multiplicity.
\par
Now, if it would be possible to develop a scattering theory in the Lax-Phillips sense for the
wave equation for $p$-forms on $H^{n+1}$, there would exist an incoming or an outgoing translation
 representation of $( H,U(t)=e^{\overline{A}t})$, that is there would exist a unitary operator
$$R : H \longrightarrow L^2({\Real},N),$$
such that
$$ R  U(t) R^{-1}= T_0(t), $$
where $N$ is an auxiliary Hilbert space and $T_0(t)$ denotes the translation of $t$ units
 to the right. In turn, this would imply that $i\overline{A}$ has a purely absolutely continuous
 spectrum on the whole real axis, and moreover that the spectrum of $i\overline{A}$ has constant
 multiplicity. As a matter of fact, if there would exist an incoming or an outgoing translation
 representation, $U(t)$ on $ H$ would be unitarily equivalent to $T_0(t)$ on $L^2({\Real},N)$,
 hence its generator $\overline{A}$ would be unitarily equivalent to the operator $\frac{d}{dx}$
 on $L^2({\Real},N)$. As a consequence, $i\overline{A}$ would be unitarily equivalent to
$i \frac{d}{dx}$, and hence, through a Fourier transform, to the multiplication operator
$$L^2({\Real},N) \longrightarrow L^2({\Real},N)$$
$$ f(x) \longmapsto x f(x),$$
whose spectrum is purely absolutely continuous, covers the whole real line and has constant
multiplicity.\par
 Since we have shown that the multiplicity of the spectrum of $i\overline{A}$ is not constant,
 we can conclude that there can be no incoming or outgoing translation representation for
$(H,U(t))$.
Hence:
\begin{thm} It is not possible to develop a scattering theory following the Lax-Phillips approach
for all $p$-forms on $H^{n+1}$.\end{thm}

\par \noindent Only in the cases $n=2p$ or $n=2p-2$ one can hope that, restricting in some way
the class of forms, it might be possible to develop a scattering theory for the wave equation
in the Lax-Phillips sense.\par
\section{Huygens' principle}
Consider the Cauchy problem (\ref{Cauchy}) for the wave equation for $p$-forms, with smooth initial data $\omega_1$, $\omega_2$, and fix a point $x_0 \in H^{n+1}$. We will say that Huygens' principle holds if the following conditions are fulfilled:
\begin{enumerate}
\item if the initial data are zero in the ball $B_{H^{n+1}}(x_0,R)$ then the solution $\omega(x,t)$ at time $t$ of the initial value problem is zero in the ball $B_{H^{n+1}}(x_0,R-t)$;\
\item if the initial data vanish outside the ball $B_{H^{n+1}}(x_0,R)$, then the solution $\omega(x,t)$ at time $t$ of the initial value problem vanishes in the ball $B_{H^{n+1}}(x_0,|t|-R)$.
\end{enumerate}
\par
\bigskip
The spectral properties of the Laplace-Beltrami operator  can be used in order to prove that Huygens' principle fails in the case of the wave equation for $p$-forms on $H^{n+1}$, when $p \not= \frac{n+1}{2}$.\par \bigskip
Consider the Cauchy problem (\ref{Cauchy}) for the wave equation for $p$-forms on $H^{n+1}$, with initial data $\omega_1, \omega_2$ smooth and compactly supported.

As noticed before, the operator $A$ defined on $\Lambda^p_c(H^{n+1}) \oplus \Lambda^p_c(H^{n+1})$ by (\ref{A}) is essentially skewselfadjoint in  $H$ endowed with the energy norm, and its closure $\overline{A}$ is the generator of a one-parameter group of unitary operators $U(t)$
on $H$ which describes the solution of the wave equation for $p$-forms on $H^{n+1}$.\par
Now, let $D_+$ be the subspace of $H$ of those smooth initial data
$ (\omega_1,\omega_2)$ such that for every $ t >0$
$$ \left\{ \begin{array}{l} (d\omega_1)(x,t)=0 \\
                   ( \delta\omega_1)(x,t)=0 \\
                    \omega_2(x,t)=0 \\ \end{array} \right.$$
for $d_{H^{n+1}}(x,x_0)<t$, and
$$(\omega_1(x,t), \omega_2(x,t))= (U(t)(\omega_1, \omega_2))(x).$$\par
The closure $\overline{D_+}$ of $D_+$ is a closed subspace of $ H$ such that
$$ U(t)\overline{D_+} \subseteq \overline{D_+}$$
for every $t>0$. To show that
\begin{equation} \label{bigcap} \bigcap_t U(t)\overline{D_+} = \{\underline{0}\} \end{equation}
note that, if $(\omega, \eta) \in U(t)\overline{D_+}$, for $t>0$, there exists a sequence
$\{(\omega_n, \eta_n)\} \subset D_+$ such that
$$U(t)(\omega_n, \eta_n) \longrightarrow (\omega, \eta)$$ in $H$, that is:
$$ \left\{ \begin{array}{l} d\omega_n \longrightarrow d\omega \\
                   \delta \omega_n \longrightarrow \delta \omega \\
                   \eta_n \longrightarrow \eta \\ \end{array} \right.$$
Now,
$$\left\{ \begin{array}{l} (d\omega_n)(x,t)=0 \\
                   (\delta \omega_n)(x,t)=0 \\
                   \eta_n(x,t)=0 \\ \end{array} \right.$$
when $d_{H^{n+1}}(x,x_0)<t$, and this implies:
$$\left\{ \begin{array}{l} (d\omega)(x,t)=0 \\
                   (\delta \omega)(x,t)=0 \\
                   \eta(x,t)=0 \\ \end{array} \right.$$
when $d_{H^{n+1}}(x,x_0)<t$.
Therefore,
if $(\omega, \eta) \in \bigcap_t U(t) \overline{D_+}$, then
$$\left\{ \begin{array}{l} (d\omega)(x,t)=0 \\
                   (\delta \omega)(x,t)=0 \\
                   \eta(x,t)=0 \\ \end{array} \right.$$
when $d_{H^{n+1}}(x,x_0)<t$ for every $t>0$, and this implies $(\omega, \eta)=\underline{0}$
in $H$, i.e. that \ref{bigcap} holds.
\par \bigskip
Now, if Huygens' principle would hold, then
$$\overline{\bigcup_t U(t)\overline{D_+}}= H.$$
In fact, if $(\omega, \eta)=0$ for $d_{H^{n+1}}(x_0,x)>R$, then, by Huygens' principle,
$$U(t)(\omega, \eta)=0 \;\; \mbox{for} \;\; d_{H^{n+1}}(x_0,x) < t-R$$
for every $t>0$.
 In particular,
$$U(t+R)(\omega, \eta)=0 \; \mbox{for} \; d_{H^{n+1}}(x_0,x)< t+R-R=t$$
for every $t>0$.
This in turn would imply
$$U(R)(\omega, \eta)\in D_+.$$
Thus, $U(-R)D_+$ would contain the set of all initial data with support in $B_{H^{n+1}}(x_0,R)$,
and hence
$$\bigcup_tU(t)D_+ \subseteq \bigcup_tU(t) \overline{D_+}$$
would contain all compactly supported data, which are dense in $H$.\par
Therefore, if Huygens' principle would hold, $D_+$ would be an outgoing subspace for $U(t)$,
and, thanks to Sinai's theorem, there would exist an outgoing translation representation for $H$,
$U(t)$. Now, this contradicts the properties of the spectrum of $\Delta_p$.
In conclusion:
\begin{thm}Huygens' principle does not hold for the wave equation for $p$-forms on $H^{n+1}$,
when $p \not= \frac{n+1}{2}$.\end{thm}

\section{Coclosed case}
In the case $p=\frac{n}{2}$, the spectrum of the Laplace-Beltrami
operator is $[0,+\infty)$, with non constant multiplicity.
This fact makes it conceivable that, under some restrictions on the
 forms considered, a scattering theory in the Lax-Phillips
sense could be developed for the wave equation.\par
This is indeed the case, if the
class of coclosed $p$-forms, with $n=2p$, in $H^{n+1}$ is
considered.\par
Indeed,  it has been shown in \cite{Branson-Olafsson-Schlichtkrull} that Huygens' principle holds for the wave equation for coclosed $p$-forms on $H^{n+1}$, when $n=2p$.

As a consequence:
\begin{thm} \label{incoming} Incoming and outgoing subspaces exist for the wave equation for coclosed $p$-forms on $H^{n+1}$, when $p=\frac{n}{2}$.
\end{thm}
Proof. Consider the space of coclosed
smooth, compactly supported $p$-forms in $H^{n+1}$, with $n=2p$, and
 the Cauchy problem (\ref{Cauchy}) for the wave equation with coclosed
smooth initial data $\omega_1$, $\omega_2$.
The condition
$$ \label{coclosed}
\delta\omega_1(x)= \delta \omega_2(x)=0 $$
is equivalent to
 $$ \forall t
\;\; \delta \omega(t,x)= \delta \omega_t (t,x)=0,$$
where
$\omega(x,t)$ is the solution of the Cauchy problem corresponding
to the initial data $\omega_1, \omega_2$.\par
Indeed,if the initial
data $\omega_1, \omega_2$ are coclosed and $\omega(x,t)$ denotes
the solution of the corresponding Cauchy problem, then $\delta \omega(x,t)$ is
the solution of the Cauchy problem for the wave equation for
$(p-1)$-forms, with initial data $\delta \omega_1, \delta
\omega_2$. Now, thanks to the uniqueness of the solution of the
Cauchy problem for the wave equation for $(p-1)$-forms, if $\delta
\omega_1 =0, \delta \omega_2=0$ then $$\delta \omega(x,t)=0$$ for
every $t$.\par
Consider now the Cauchy problem for the wave equation in the space
of coclosed forms. The energy left invariant by the system is expressed
by the norm
$$ \|
(\omega, \eta)\|^2_{\tilde{H}}= < \Delta_p \omega, \omega> + \|\eta\|^2= \|d
\omega\|^2 + \| \eta\|^2  \geq 0, $$
which is induced by the scalar product (\ref{scalarproduct}) restricted to the set of
coclosed forms.
The Hilbert space $\tilde{{H}}$, completion of the space $ D$ of all smooth, compactly
supported coclosed data
with respect to this norm, is a closed subspace of ${H}$. \par

We will denote by  $\tilde A$ the restriction of the operator $A$, expressed by (\ref{A}),
to $ D$.  $\tilde A$ maps $ D$ into $ D$, and is essentially skewselfadjoint in $\tilde{ H}$.
 Consequently, it is the generator of a one-parameter group of unitary operators $\tilde{U}(t)$,
 which describes the evolution of the system,
 mapping initial data to the corresponding solution at time $t$.

\bigskip

Let now $\tilde{D}_+$ be the subspace of $\tilde{{H}}$ of those initial data
$(\omega_1, \omega_2)$ such that for every $ t >0$
$$\left\{ \begin{array}{l} (d\omega_1)(x,t)=0 \\
                   \omega_2(x,t)=0 \\ \end{array} \right.$$
for $d_{H^{n+1}}(x,x_0) <t$, where $x_0$ is a fixed point of $H^{n+1}$ and
$$(\omega_1(x,t), \omega_2(x,t))=(\tilde{U}(t) (\omega_1, \omega_2))(x).$$
Its closure $\overline{\tilde{D}}_+$ in $\tilde{H}$ is a closed subspace of $\tilde{ H}$ such that
$$\tilde{U}(t)\overline{\tilde{D}_+} \subseteq \overline{\tilde{D}_+}.$$
Moreover, it satisfies
$$ \bigcap_t \tilde{U}(t)\overline{\tilde{D}_+} = \{\underline{0}\}.$$

Finally, the validity of Huygens' principle implies that
$$\overline{\bigcup_t \tilde{U}(t)\overline{\tilde{D}_+}}= \tilde{ H},$$
because, if Huygens'principle holds, then, for every $R>0$, $\tilde{U}(-R)\tilde{D}_+$ contains every smooth compactly supported coclosed data with support in $B(x_0,R)$.
Hence,
$\overline{\tilde{D}_+}$ is an outgoing subspace for $\tilde{U}(t)$ in $\tilde{ H}$.
Analogously, the subspace $\overline{\tilde{D}_-}$ of $\tilde{{H}}$, consisting of those initial data $(\omega_1, \omega_2)$ such that for every $ t <0$
$$\left\{ \begin{array}{l} (d\omega_1)(x,t)=0 \\
                   \omega_2(x,t)=0 \\ \end{array} \right.$$
for $d_{H^{n+1}}(x,x_0) <-t$, can be shown to be an incoming subspace for $\tilde{U}(t)$ in
$\tilde{ H}$.\qed
\par \bigskip
Theorem \ref{incoming} implies that the spectrum of $i\overline{\tilde A}$ is purely absolutely continuous, with constant multiplicity, and equal to $(- \infty, +\infty)$. As a consequence:
\begin{thm}
When $n=2p$, the spectrum of the Laplace-Beltrami operator on coclosed $p$-forms on $H^{n+1}$ is purely absolutely continuous, with constant multiplicity, equal to $[0,+\infty)$. \end{thm}

\end{document}